# A Lane Usage Strategy for General Traffic Access on Bus Lanes under Mixed Traffic Environment

Haoran Li, Graduate Student Member, IEEE, Zhenzhou Yuan, Rui Yue, Guangchuan Yang, Chuang Zhu, Graduate Student Member, IEEE, and Siyuan Chen

*Abstract*—The strategy of permitting general traffic to use the bus lane for improved utilization while ensuring bus priority has gained increasingly attention, particularly with the support of vehicle-to-everything technology. In this study, we propose a novel lane usage strategy called Dynamic Spatial-Temporal Priority (DSTP) to ensure bus priority and optimize bus lane usage in a mixed traffic environment. DSTP leverages dynamic methods to identify available spatial-temporal resources in the lane, utilizing signal timing, road information, and vehicle data. A Right-of-Way assignment optimization model is then developed based on these resources to determine which vehicles can enter the bus lane. The model is dynamically enacted using a rolling horizon scheme to accommodate time-varying traffic conditions. Numerical studies have validated the advantages of DSTP, showing maintained bus priority, improved traffic efficiency, reduced fuel consumption, and lower $CO_2$ emissions, especially during periods of high traffic demand and concentrated bus arrivals.

## I. INTRODUCTION

REASONABLE public transport priority (PTP) measures have the potential of improving the travel time and reliability of public transportation systems [1]. One of the primary types of PTP measures is providing road space priority, which is typically realized through the implementation of bus lanes. Bus lanes enhance the efficiency of public transportation by minimizing traffic interactions between buses and private vehicles, ensuring stable and reliable service for commuters [2] and contributing to environmental sustainability by reducing carbon emissions and other pollutants emitted by buses [3].

Regrettably, constrained urban road resources pose a challenge as converting general lanes into bus lanes may exacerbate traffic congestion by diminishing the overall capacity for general traffic [4]. Bai et al. [5] argued that while the bus lane reservation strategy effectively reduces carbon emissions and enhances bus service quality at high utilization rates, it may have the opposite effect at low utilization rates.

This research was supported by the Fundamental Research Funds for the Central Universities under Grant 2023YJS038, in part by the Beijing Natural Science Foundation under Grant J210001. (Corresponding authors: Rui Yue; Chuang Zhu.)

Haoran Li, Zhenzhou Yuan, Rui Yue, and Chuang Zhu are with the Key Laboratory of Transport Industry of Big Data Application Technologies for Comprehensive Transport, Beijing Jiaotong University, Beijing 100044, China (e-mail: lihaoran@bjtu.edu.cn; zzyuan@bjtu.edu.cn; yuerui@bjtu.edu.cn; zhuchuang@bjtu.edu.cn).

Guangchuan Yang is with the Institute of Transportation Research and Education, North Carolina State University, Raleigh, NC 27606, USA (e-mail: gyang24@ncsu.edu)

Siyuan Chen is with the Zhejiang Communications Investment Group Co., Ltd., Hangzhou, Zhejiang 310000, China (e-mail: chensienna@163.com)

These issues can be addressed by directing general traffic to intermittently utilize the bus lane during periods of bus absence, coupled with some rules to ensure priority for buses upon arrival. Viegas and Lu [6, 7] introduced the concept of the intermittent bus lane (IBL) to allow general traffic to enter the bus lane when there are no buses present. This strategy requires coordination with the modification of traffic signals. Eichler improved upon the concept of IBL and introduced the Bus Lane with Intermittent Priority (BLIP) [8]. BLIP prohibits general traffic from entering the bus lane when a bus enters the road, and requires general traffic preceding the bus to vacate the bus lane upon the approach of buses. This strategy does not require changes to traffic signals, resulting in minimal impact on the other approach of the intersection. A considerable amount of research has utilized theoretical models [9, 10] or simulation-based approaches [11-13] to evaluate the applicability and transferability of the BLIP strategy. Initially, the BLIP strategy considered clearing all vehicles ahead of a bus. To address this issue, Ma et al. [14] proposed dividing road sections into multiple empty areas using Variable Message Signs (VMS). When a bus enters the corresponding empty area, general traffic within that area is required to exit the bus lane.

Implementing IBL/BLIP requires conveying information to drivers through VMS or vertical pole lights, which limits the effectiveness of these strategies in controlled roads [15, 16]. With the development of V2X, real-time lane information can be transmitted to drivers' mobile navigation or On-Board Equipment (OBE) via various communication methods such as dedicated short-range communication (DSRC), cellular network, or hybrid [17, 18]. Levin et al. [19] suggested that V2X technologies could be used to quickly inform general traffic whether they can use the bus lane, but in their strategy, all general traffic ahead of a bus still needs to exit. Wu et al. [16] further enhanced the BLIP strategy by incorporating V2X technologies, introducing the Bus Lane Intermittent and Dynamic Priority (BLIDP): as a bus travels along the road, it broadcasts V2V information, and general traffic within a specific clearance distance ahead of the bus must exit the bus lane. Ou et al. [20] also adopted a similar methodology for application in tram lanes. Luo et al. [21] introduced dynamic bus lane with moving block in a connected environment, which can adjust the clearance distance based on the bus speed. Othman et al. [15] compared the BLIDP strategy with exclusive bus lanes and mixed lanes (without bus priority) across varying traffic demand and bus frequencies, assuming all cars are connected vehicles. The results indicated that the BLIDP strategy demonstrated superior performance under intermediate levels of traffic demand. However, Othman's findings also revealed that under high traffic demand levels, there was little difference in the average bus speed between the BLIDP strategy and the mixed lanes strategy, highlighting the limitations of BLIDP in such scenarios. The phenomenon was previously explored by Kampouri [22], who observed that the BLIDP strategy can exacerbate congestion within the bus lane during periods of high traffic demand. Specifically, general traffic that is required to vacate the bus lane may be forced to block bus priority due to the difficulty in finding timely lane-changing opportunities in adjacent high-traffic lanes. Some studies suggest that improving bus lane utilization could also involve allowing a portion of general traffic to enter the bus lane, creating a controlled mixed traffic lane. Chen et al. [23] proposed the transformation of Bus Rapid Transit (BRT) lanes into mixed-use lanes, allowing AVs to utilize the BRT lanes in the early stage of AV adoption. Anderson et al. [24] proposed a "dynamic bus lane strategy" where the control system adjusts the proportion of vehicles to buses in one lane of a multi-lane arterial road instead of enforcing complete separation of vehicle types. Shao et al. [25] considered a strategy for controlling connected and human-driven vehicles (CHVs) actively borrowing dedicated connected and automated vehicle (CAV) lanes on highways. Although this study did not specifically focus on bus lanes, they thought this strategy could also apply to bus lanes. Shan et al. [26] proposed a trajectory optimization method that allows some CAVs to enter the bus lane and optimize their own trajectories based on the known bus trajectories.

The discussion on improving bus lane utilization has never ceased, and previous studies have made significant contributions to this field. However, existing studies still exhibit notable limitations. Regarding strategies that require clearing the general traffic in front of buses, while they theoretically ensure absolute priority for buses, as mentioned earlier, with the increasing level of traffic demand, the actual priority of buses cannot be guaranteed due to the inability to promptly clear general traffic on the bus lane. Regarding strategies that allow buses to share lanes with some general traffic, the absolute priority of buses also cannot be guaranteed. This is because, at signalized intersections, buses may be delayed from passing through the intersection within the green light time due to obstacles from the general traffic ahead. Therefore, the first limitation of existing research is the lack of strategies that can practically ensure the absolute priority of buses, especially under high traffic demand levels. Furthermore, the potential weaving issues on bus lanes deserve attention. Past studies primarily focus on accommodating through traffic within bus lanes and do not adequately consider the needs of right-turn traffic. However, right-turn traffic on general lanes often intersects with through traffic on bus lanes because right-turn traffic typically needs to use the bus lane to make a right turn or cross the bus lane to enter a right-turn lane [27]. The weaving issues become more pronounced as the general traffic using bus lanes increases. There is limited research testing whether strategies proposed are still effective in scenarios with significant right-turn traffic. Researchers have not extensively studied how strategies perform when faced with a lot of right-turn traffic. Moreover, Past research has often overlooked bus stops as a factor in strategy development. For instance, in strategies that involve clearing general traffic, when a bus is at a stop, even if it has a long dwell time, general traffic is not allowed to enter the clearance distance ahead, leading to prolonged waste of road resources. Consequently, it is crucial to construct a strategy that considers a more comprehensive set of factors.

With these concerns, this paper introduces a nuanced strategy named dynamic spatial-temporal priority (DSTP), aiming at filling the aforementioned gaps. The main idea is that a few vehicles are allowed to enter bus lanes when there are available spatial-temporal resources under the cooperative vehicle infrastructure system. The spatial-temporal resources refer to the remaining road space resources in the bus lane under the condition of ensuring priority passage for buses, along with the green light time resources at intersections, which will be identified based on the method we propose while considering the impacts of bus stops, signal timings, and vehicle states. A Right-of-Way assignment optimization model is constructed to match the optimal CHVs or CAVs utilizing these resources. The goal of proposed strategy is to ensure that the bus lane is utilized in a manner that minimally interferes with bus operations. The main contributions of this study are summarized as follows: (1) Developing a novel strategy to enhance bus lane utilization called DSTP under a partially connected vehicle environment. DSTP incorporates multiple modules and operational rules that permit general traffic to enter the bus lane without a mandatory exit. Importantly, this strategy considers protocols for both through traffic and right-turn traffic to access the bus lane. (2) Proposing a method to dynamically identify the remaining spatial-temporal resources in the bus lane while ensuring the normal operation of buses, considering the collective influence of bus stops, signal timings, and vehicle conditions. (3) Constructing a Right-of-Way assignment optimization model to determine the optimal vehicles capable of utilizing idle bus lane resources. (4) Analyzing the benefits of the novel DSTP strategy and other strategies across different road and traffic flow parameters.

The organizational structure of this paper is described as follows. Section II outlines the problems to be addressed by this study. Section III presents the proposed strategy in detail. Section IV describes the simulation experiments conducted, and evaluates the performance of the proposed DSTP strategy through a comparative analysis with the Exclusive Bus Lane (EBL) and the BLIDP strategies. Section V concludes the paper and discusses future research directions.

## II. PROBLEM DESCRIPTION

Fig. 1 (a) illustrates the signalized roadway intersection with some roadside units (RSUs) under connected environment in this study. There exist two types of road lanes, i.e., general lane for general vehicles and no-barrier-separated bus lane for buses, and consideration is limited to traffic processing in one direction (southbound). There is a no-changing zone located in the vicinity the stop bar, prohibits lane change activities within its confines. As a common practice, vehicles must complete their lane changes before entering the no-changing zone. CAVs, HDVs and CHVs coexist in the general lanes and follow the intersection signals. The control center disseminates lane-changing advisories to CAVs and CHVs via information transmission. In this study, the notation $(i,j)$ is used to designate vehicles within the research area, where $j$ denotes the $j$-th lane and $i$ denotes the $i$-th vehicle in the lane. The vehicles are numbered sequentially according to the time they are approaching the intersection. The lane indices $j$ are arranged in ascending order from left to right. $I_j$ is the set of vehicles in a lane $j$. The set of $j$ is $J$. Vehicles within the controlled zone are assumed to be connected. The location and speed of the $(i,j)$ are denoted by $l_{(i,j)}$ and $v_{(i,j)}$, respectively.

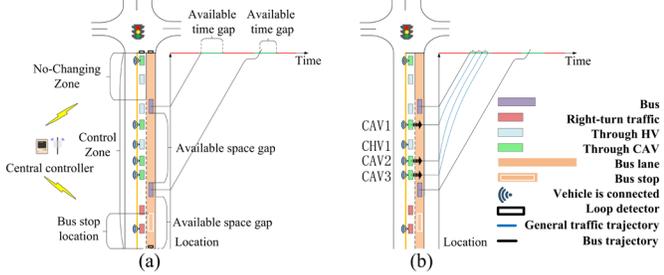

Fig. 1. Concept of DSTP: (a) details of the road; (b) identifying remaining spatial-temporal resources and conducting Right-of-Way optimization.

The expected bus operation is as depicted in Fig. 1(a), the first bus passed through the intersection without stopping during the green time, while the second bus arrived at the stop bar and came to a halt during the red time. At the present moment, it is evident that there is available space in the bus lane as well as green light time. Our goal is to dynamically identify these remaining resources on the bus lane in real-time and enable vehicles in the general lane to make reasonable use of these resources. With the implementation of the proposed strategy, the control center utilizes V2X devices to obtain vehicle status information and signal timing data. The departure time at the stop bar will be estimated for vehicles in the bus lane and the general lanes. To identify the remaining resources, we propose the concept of Available Spatial Gap (ASG), and Available Temporal Gap (ATG) to improve the utilization of the bus lane from both spatial and temporal perspectives, where $(i_f, j)$ denotes the following vehicle of vehicle $i$ in lane $j$. An ASG is defined as the spatial gap between adjacent vehicles in the bus lane at a time step where there is enough space to accommodate at least one vehicle length plus the minimum space headway with target leader and target follower as shown in Fig. 1(a). There can be multiple ASGs in the bus lane. The ATGs are defined as the temporal gap that intersects with the green time duration and the departure times of adjacent vehicles, which can be considered as companions of ASGs.

After identifying the ASGs and their corresponding ATGs, the control center will recognize the vehicles that satisfy the ASGs and ATGs and consider them as the alternative sets of vehicles, such as CAV1, CAV2, CAV3, and CHV1 depicted in Fig. 1(b). We suppose that CAVs use Cooperative Adaptive Cruise Control (CACC) with a time headway of $\tau_c$ when following other CAVs [28]. If a CAV follows an HDV or CHV, the car-following model of the CAV will automatically switch from CACC to Adaptive Cruise Control (ACC) [29] with a time headway of $\tau_a$. Additionally, due to the significant uncertainties associated with HDVs and CHVs being manually driven, we assume that they maintain a larger time headway when following other vehicles, denoted as $\tau_h$. Then, based on the proposed optimization model that takes into account vehicle departure time, current speed, and position, the control center will attempt to find vehicle platoons that can most efficiently utilize the remaining resources and assign the Right-of-Way to these vehicles, as shown in Fig. 1(b). After determining the optimal vehicles to utilize the bus lane, the control center sends lane-changing advisories. The drivers receive these advisories in both text and audio formats, which is similar to many real-world Cooperative Intelligent Transport Systems (C-ITS) applications[30]. Meanwhile, we consider the applicability of the proposed strategy to right-turn traffic. Right-turn traffic is given higher priority to enter the bus lane, and the control center implements measures to promote the likelihood of successful lane changes. The specific protocols for implementing these measures will be described in detail in Section III.

In order to effectively implement this concept, a comprehensive and adaptable control framework, coupled with an optimization approach, must be devised. These components should seamlessly integrate with real-time roadway data, including the dimensions of the control and no-changing zones, bus stop locations, as well as signal and vehicle status information.

## III. METHODOLOGY

The general framework of the proposed strategy is illustrated in Fig. 2. The initialization module initiates the process by gathering road, signal, and vehicle status information, enabling the control center estimating the departure time of vehicles. Subsequently, this information is transmitted to the recognition module for the initial identification of ASGs. Then, the control center applies the lane-changing protocol specific to right-turn traffic to filter out eligible vehicles and save their indices for further reference. If an ASG exists, a further identification of ATG is performed. When ATG exists and the time aligns with the preset control interval, $h$, the ASG and ATG information are transmitted to the optimization module. In the optimization module, vehicles that meet the lane-changing protocol for through traffic are identified, and the optimization model is executed to select the most suitable lane-changing vehicles, with their indices being saved for later use. Finally, the saved indices of right-turn traffic and through traffic are transmitted to the execution module, and the control center issues lane-changing advisories to these vehicles.

### A. Assumption

To facilitate the model development process, several foundational assumptions are posited as follows:
- Each vehicle is willing to comply with the lane-changing instructions sent by the control center unless the driver perceives them as unsuitable given the prevailing driving conditions.
- The desired speed of general traffic equals the speed limit of the road segment. These vehicles travel at the speed limit unless impeded by the presence of a preceding vehicle.
- Communication-related issues, such as delays and data packet loss, are ignored.
- Loop detectors are installed within the communication range.

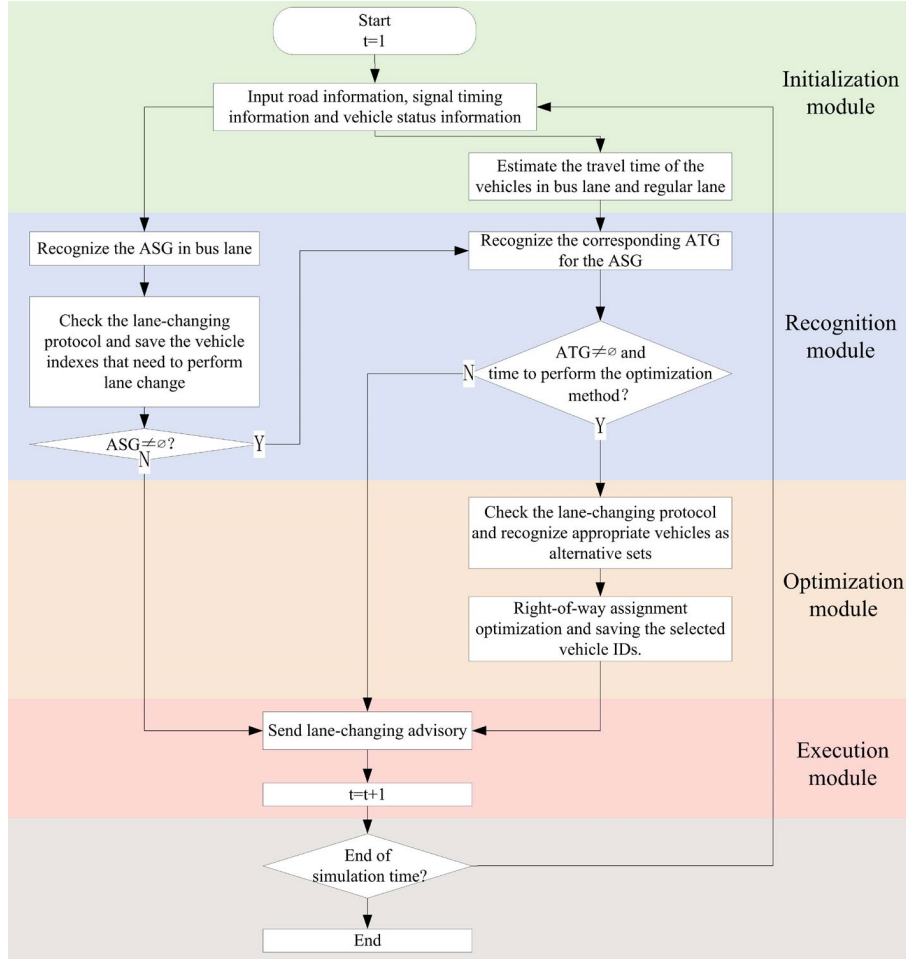

Fig. 2. Flow chart of the DSTP strategy

- The general traffic consists of vehicles with uniform size and vehicle dynamics, including acceleration and desired speed. All vehicles are light vehicles with no time delay in their acceleration and deceleration capabilities. Additionally, all following vehicles adhere to specific car-following rules, where vehicles with larger gaps proactively seek to converge with their leading vehicles while maintaining safety constraints.

*B. Vehicle departure time estimation*

To provide a basis for the implementation of the recognition and optimization modules, the departure time of vehicles within each lane is estimated individually using the kinematics method. The estimation process adheres to the subsequent procedural steps.

1) The timestamp $t_{(pre,j)}$ will be updated in real-time and recorded by the control center when the loop detectors at the stop bar on the lane j detect a vehicle passing through.

2) Defines the set of vehicles $I'_j$ that need to be estimated departure time. If lane j is a bus lane with a bus stop, and there is a bus $(\partial, j)$ between the bus stop and the entrance to the control zone, the set of vehicles to be considered for estimation should include the vehicles between $(\partial, j)$ and the stop bar, i.e., $I'_j = \{(i,j)|p_{(\partial,j)} \leq p_{(i,j)} \leq l_c, i \in I_j\}$. Otherwise, the required set of vehicles to be estimated is all the vehicles on the lane, i.e., $I'_j = I_j$.

3) Calculate the departure time $t^d_{(i,j)}$ for each vehicle to pass through the stop bar. Firstly, the departure time of the vehicle without considering the limitation of signal lights, $t^p_{(i,j)}$, is calculated as follows [31, 32]:

$$t^p_{(i,j)} = \begin{cases} max(t_{(pre,j)} + \tau, t^f_{(i,j)}), if\ i\ is\ the\ first\ vehicle \\ max(t^d_{(i_p,j)} + \tau, t^f_{(i,j)}), otherwise \end{cases}$$
(1)

$t^f_{(i,j)}$ is the earliest time of vehicle drives pass the stop bar without considering its preceding vehicle and signal control. According to the kinematic method:

$$t^f_{(i,j)} = t_0 + \frac{l_c - x_{(i,j)} - \left(\frac{(v_d^2 - v_{(i,j)}^2)}{2a_U}\right)}{v_d} + \frac{v_d - v_{(i,j)}}{a_U}$$
(2)

$\tau$ is the desire time headway determined by the type of vehicle $(i,j)$:

$$\tau = \begin{cases} \tau_a + \varepsilon_a, when\ i \in I_j^a, i_p \in I_j^h \\ \tau_c + \varepsilon_c, when\ i \in I_j^a, i_p \in I_j^a \\ \tau_h + \varepsilon_h, when\ i \in I_j^h \end{cases}$$
(3)

$\varepsilon_a$, $\varepsilon_c$, and $\varepsilon_h$ is the redundant time that considering the inaccurate driver behavior modeling. Once the $t^p_{(i,j)}$ of a vehicle is obtained, the impact of signal timing on the vehicle's movement is considered. Therefore, the departure time $t_{(i,j)}^d$ of the vehicle $(i,j)$ can be calculated as follows.

$$t_{(i,j)}^d = max(t_{(i,j)}^p, t_{(i,j)}^G + t_l) \quad (4)$$
$$t_{(i,j)}^G = t_{(i,j)}^R + t_r \quad (5)$$
$$t_{(i,j)}^R = \left\lceil \frac{t_{(i,j)}^p}{t_c} \right\rceil \cdot t_c \quad (6)$$
$$t_c = t_r + t_g \quad (7)$$

where $t_l$ is the start-up lost time, which depends on the vehicle's position in the queue. $t_{(i,j)}^R$ is the start of the red time of the signal cycle in which vehicle $(i,j)$ arrives at the stop bar. $t_{(i,j)}^G$ is the start of the green time of the signal cycle in which vehicle $(i,j)$ arrives at the stop bar. It should be noted that when lane $j$ is a through and right turn lane while vehicle $(i,j)$ is a right-turn vehicle, its departure time is not affected by signal light restrictions, i.e., $t_{(i,j)}^d = t_{(i,j)}^p$.

Specifically, when vehicle $(i,j)$ is a bus stopped at a bus stop, the calculation method for its departure time is as (8) and (9) at the bottom of the page, where $t_s$ represents the dwell time:

$$t_{(i,j)}^d = max(t_{(i,j)}^p, \left\lceil \frac{t_{(i,j)}^p}{t_c} \right\rceil \cdot t_c + t_r + t_l) \quad (8)$$

*C. Recognition module*

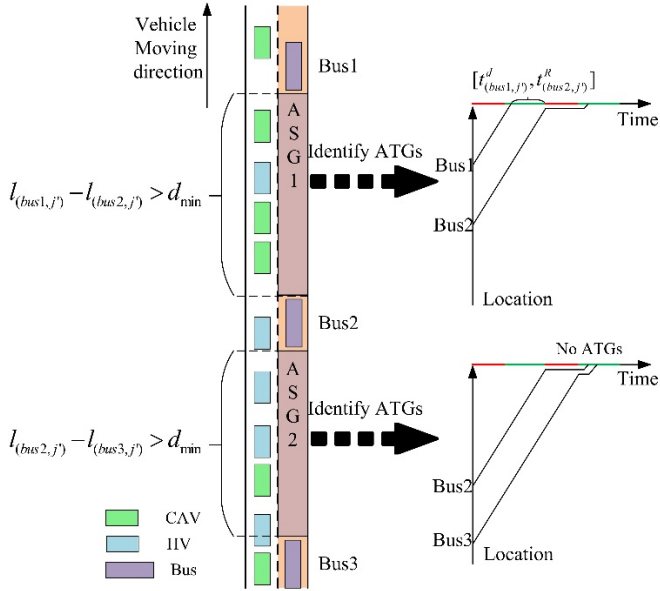

Fig. 3. Recognition process of ASGs and ATGs

The recognition module is utilized to recognize ASGs and ATGs of the bus lane as shown in Fig. 3. Firstly, the ASGs of the bus lane is determined by calculating the gap between the preceding and following vehicles in order in the bus lane. Next, for each ASG, the ATGs are identified based on the departure time of vehicles. ATGs can be considered as usable green time windows that can be utilized without disrupting the vehicles passing through the stop bar in the bus lane.

At each time step, the control center calculates the spatial gaps between vehicles within the control zone according to their proximity to the intersection. Initially, the range recognized by ASG is initialized as $[0, l_c]$. When a bus is approaching or has already stopped at the bus stop, the starting point for ASG recognition shifts to $[p_{(\partial,j)}, l_c]$.

For longitudinal safety concerns, if the defined ASG is too small, collisions can occur between vehicles in adjacent lanes and the leading or following vehicles of the ASG during lane changes. Therefore, the minimum usable length of an ASG is calculated as follows, which must have enough space to accommodate at least one vehicle for a lane change.

$$d_{min} = d_p + l_v + d_f \quad (10)$$

$d_{min}$ represents the minimum gap required for an ASG. $l_v$ represents the length of a vehicle. $d_p$ and $d_f$ represent the safe distance between the new lead vehicle or the following vehicle changing lanes onto the bus lane. In fact, the specific values of $d_p$ and $d_f$ depend on the vehicles associated with the ASG and the driving conditions of the vehicles intending to change lanes [33, 34]. In this module, we set $d_p$ and $d_f$ as fixed values for simplification. The exact modeling of longitudinal safety with lane changing is presented in the optimization module in Section III, E. To handle the inaccurate modeling of these safety concerns, a rolling horizon scheme is proposed in Section III, F.

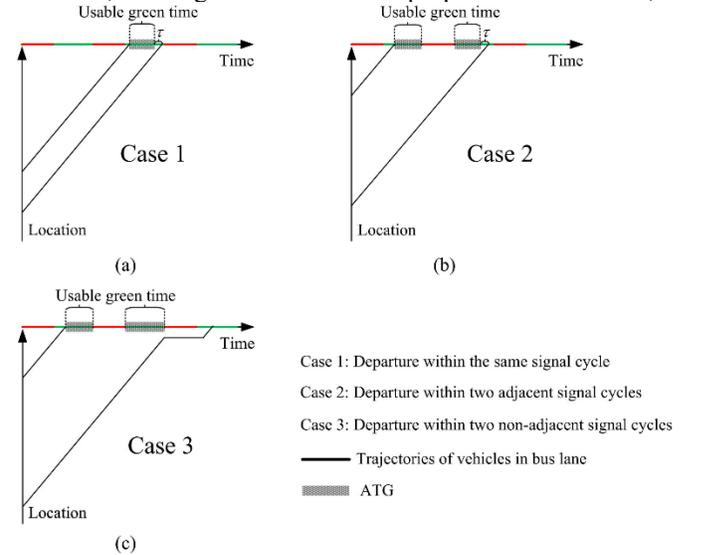

Fig. 4. The cases exist when recognizing ATGs

If $l_{(i,j)} - l_{(i_f,j)} \geq d_{min}$, it indicates that an ASG has been recognized, and its corresponding ATG will be calculated. ATG is defined as the remaining green signal time between the departure time $t_{(i,j)}^d$ of the preceding vehicle, $(i,j)$, and the departure time $t_{i_f,j}^d$ of the following vehicle, $(i_f, j)$. To mitigate the impact on the departure time $t_{(i_f,j)}^d$ of vehicle $(i_f, j)$, the end of the time window has been advanced and reduced by $\tau$ seconds, i.e., $[t_{(i,j)}^d, t_{(i_f,j)}^d - \tau] \cap \xi$. Where $\xi$ is the collection of

$$t_{(i,j)}^p = max(t_{(i_p,j)}^d + \tau, t + \frac{(l_c - l_s) - \left(\frac{(v_d)^2 - v_{(i,j)}^2}{a_U}\right)}{v_d} + \frac{v_d - v_{(i,j)}}{a_U} + t_s) \quad (9)$$

green phases. Due to the possibility of the departure times of the two vehicles not falling within the same signal cycle, the ATG may consist of multiple time windows, as shown in Fig. 4. To facilitate the description of each time window within the ATG, we designate $T_n(n = 1,2,\ldots,N)$ as the nth time window in the ATG. The algorithm for recognizing ATG is as in Algorithm 1.

**Algorithm 1** ATG reorganization on bus lane
---
Input: $t$, $t^d_{(i,j)}$, $t^d_{(i_f,j)}$, $\tau$, $t_r$, $t_c$.
Output: $T_n$

Initialization: $N \leftarrow \left\lfloor \frac{t^d_{(i_f,j)}}{t_c} \right\rfloor - \left\lfloor \frac{t^d_{(i,j)}}{t_c} \right\rfloor + 1$ , $T_n \leftarrow \emptyset (n = 1,2,\ldots,N)$
If $N = 0$ & $t^d_{(i_f,j)} - \tau - t^d_{(i,j)} > \tau_a$
$\quad T_1 \leftarrow [t^d_{(i,j)}, t^d_{(i_f,j)}]$
Else if $N = 1$
If $t^d_{(i,j)} - (\left\lfloor \frac{t^d_{(i,j)}}{t_c} \right\rfloor + 1) \cdot t_c \geq \tau_a$
$\quad T_1 \leftarrow [t^d_{(i,j)}, (\left\lfloor \frac{t^d_{(i,j)}}{t_c} \right\rfloor + 1) \cdot t_c]$
If $t^d_{(i_f,j)} - \tau - (\left\lfloor \frac{t^d_{(i_f,j)}}{t_c} \right\rfloor \cdot t_c + t_r) \geq \tau_a$
$\quad T_2 \leftarrow [(\left\lfloor \frac{t^d_{(i_f,j)}}{C} \right\rfloor \cdot t_c + t_r), t^d_{(i_f,j)} - \tau]$
Else if $N > 1$
$\quad T_n \leftarrow [(\left\lfloor \frac{t^d_{(i,j)}}{t_c} \right\rfloor + n) \cdot t_c + t_r, (\left\lfloor \frac{t^d_{(i,j)}}{t_c} \right\rfloor + n + 1) \cdot t_c], n = 2,3,\ldots N-1$
If $t^d_{(i,j)} - (\left\lfloor \frac{t^d_{(i,j)}}{t_c} \right\rfloor + 1) \cdot t_c \geq \tau_a$
$\quad T_1 \leftarrow [t^d_{(i,j)}, (\left\lfloor \frac{t^d_{(i,j)}}{t_c} \right\rfloor + 1) \cdot t_c]$
If $t^d_{(i+1,j)} - \tau - (\left\lfloor \frac{t^d_{(i_f,j)}}{t_c} \right\rfloor \cdot t_c + t_r) \geq \tau_a$
$\quad T_N \leftarrow [(\left\lfloor \frac{t^d_{(i_f,j)}}{t_c} \right\rfloor \cdot t_c + t_r), t^d_{(i_f,j)} - \tau]$

*D. Lane-changing protocol for right-turn and through traffic*

This section is intended to describe the lane-changing rules, considering the scenario where right-turn vehicles are allowed to use the bus lane for their right turns while also considering whether through traffic can utilize the bus lane. The specific lane-changing protocol is outlined as follows: The protocol of the lane-changing is constructed mainly to: 1. Ensure that right-turn vehicles can utilize the bus lane for their right turns. 2. After ensuring the proper accommodation of right-turn vehicles, consider whether through traffic can use the bus lane. To facilitate clear distinction, the bus lane is designated as $j'$. The specific lane-changing protocol is outlined as follows:

1) For each ASG composed of preceding vehicle $(p, j')$ and following vehicle $(f, j')$, when there are connected right-turn vehicles on the adjacent general lane within their range, they will be advised to change lanes. The corresponding set of vehicles is denoted as $I_j^{advice\_1} = \{(i,j) | l_{(f,j')} < l_{(i,j)} < l_{(p,j')}, i \in I_j^r\}$, where $I_j^r$ is the right-turn vehicles in lane $j$.

2) For those connected right-turn vehicles that have reached a distance $l_r$ from the stop bar but have not yet received a lane-changing advisory, the control center will notify them to autonomously choose an appropriate moment to change lanes. The corresponding set of vehicles is denoted as $I_j^{advice\_2} = \{(i,j) | l_{(i,j)} > l_c - l_r, i \in I_j^r\}$. Therefore, the set of right-turn vehicles that receive the lane-changing advisory is $I_j^{advice\_1} \cup I_j^{advice\_2}$.

3) Grant the highest priority of using the bus lane to non-connected right-turn vehicles. If a bus stop is present on the road, they are required to enter the bus lane immediately after passing the bus stop; otherwise, they are required to enter the bus lane as soon as they enter the road entrance.

4) For each ASG, when there are no connected right-turn vehicles ahead in the adjacent lane of the preceding vehicle, i.e., when $\{(i,j) | l_{(f,j')} < l_{(i,j)} < l_c, i \in I_j^r\} = \emptyset$, indicates that the through traffic within the ASG range have satisfied the lane-changing protocol. When the ASG has a corresponding ATG, the control center begins to recognize appropriate vehicles as alternative sets. To avoid stop-and-go traffic in the bus lane, the control center includes vehicles in the alternative sets that are positioned within the ASG range and can pass through the intersection within the available time window in the adjacent lanes, i.e., $I_j^X = \{(i,j) | l_{(f,j')} \leq l_{(i,j)} \leq l_{(p,j')}, t^f_{(i,j)} \in T_n, i \in I_j\}$.

*E. Optimization module*

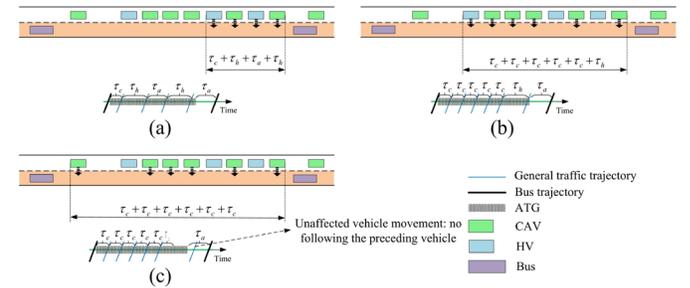

Fig. 5. Results of the Right-of-Way assignment decision process: (a) without optimization model implementation; (b) with optimization model implementation; (c) the last CAV cannot closely follow the preceding vehicle.

In this section, a Right-of-Way assignment optimization model is introduced to determine which vehicles in the set $I_j^X$ are eligible to enter the bus lane. Specifically, given that the green time is fixed, following a sequential order of vehicles without introducing an optimization method would result in only the first four vehicles fully utilizing the available green time, as shown in Fig. 5(a). However, by integrating an optimization approach, a more efficient solution based on the following characteristics of the vehicles as described in Section II is derived, as shown in Fig. 5(b). Additionally, when constructing the proposed optimization model, it does not only take into account the vehicles' following characteristics, as not every vehicle passes the stop bar in a following manner, as demonstrated by the last CAV in Fig. 5(c).

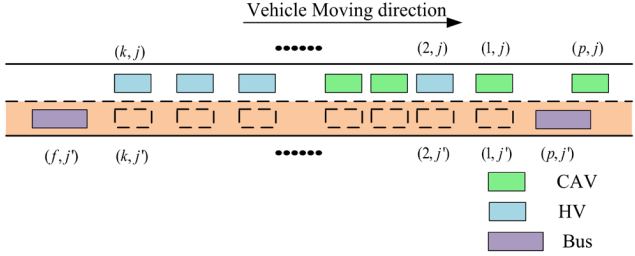

Fig. 6. Anchoring vehicle numbering using virtual vehicles

Once $I_j^X$ is determined, the control center allocates the Right-of-Way to suitable general traffic within ASG to reduce the total travel time in the $I_j^X$. To facilitate model construction and solving, we continue to consider the concept of virtual vehicles to anchor the vehicle indices temporarily. This consideration is primarily based on the following: Since vehicle indices are named based on the order in the lane, lane changes by vehicles can result in changes to their indices. Once the optimization model enters the iterative process, these changes in vehicle indices can render the optimization model ineffective. Therefore, we renumber the vehicles in the order of increasing distance from the stop bar within $I_j^X$, creating a new index set $K = \{k | k = 1, 2, \ldots, |I_j^X|\}$. Then, as shown in Fig. 6, suppose that there is an ASG and ATG between vehicles $(p, j')$ and $(f, j')$ in the bus lane. Vehicle $(k, j)$ is considered to have a corresponding virtual vehicle $(k, j')$ in the bus lane $j'$. If $(k, j)$ is chosen to change lanes, $(k, j')$ becomes a real vehicle, and $(k, j)$ becomes a virtual vehicle. In this approach, the indices of vehicles on the lane will remain unchanged even when they change lanes. The virtual vehicles inherit the state (position, velocity, departure time) of the nearest preceding real vehicle to estimate the departure time of the subsequent real vehicles. Although there is a time delay for vehicles to change lanes and occupy the corresponding virtual space, it does not affect our estimation of their longitudinal travel time. This is because the lateral and longitudinal movements of vehicles are synchronized, and we consider the influence of the preceding vehicle's operational state on the departure time of the vehicle in the constraints. Additionally, we incorporate lateral safety constraints to ensure safe lane changes.

In this study, we aim to minimize the travel time of vehicles in the objective function, achieving the allocation of Right-of-Way to general traffic, allowing them to utilize the bus lane until they pass the intersection:

$$min \sum_{k \in K}(1 - x_k) \cdot t_{(k,j)}^d + \sum_{k \in K} x_k \cdot t_{(k,j')}^d \quad (11)$$

Where $x_k$ is a binary variable, and if $x_k = 1$, it represents sending a lane-changing advisory to vehicle $(k, j)$.

1) Longitudinal vehicle operation constrains

For vehicles in general lane $j$, if vehicle $(k, j)$ changes lane, $(k, j)$ becomes a virtual vehicle with a departure time equal to the preceding vehicle. Otherwise, it is an actual vehicle whose departure time is the value estimated in the initialization module. Therefore, as shown in (12) at the bottom of the page, for the first vehicle $(1, j)$, if it does not change lanes, its departure time is the maximum of the free-flow travel time (i.e., $t_{(1,j)}^f$) and the travel time following the preceding vehicle (i.e., $t_{(p,j)}^d + \tau$); whereas if it changes lanes, the departure time of vehicle $(1, j)$ is the same as that of the preceding vehicle. $\eta_1$ is a binary auxiliary variable used to determine whether $(1, j)$ is also the first vehicle on lane $j$. If that is the case, there are no vehicles ahead of vehicle $(1, j)$ and it can pass the stop bar with free-flow travel time, as shown in (13).

$$\eta_1 = \begin{cases} 1, if \ (t_{(p,j)}^d + t_g) mod(t_g + t_r) \neq 0 \\ 0, others \end{cases} \quad (13)$$

The departure time of the vehicles behind vehicle $(1, j)$ as indicated in (14). The departure time of vehicle $(k, j')$ in bus lane $j'$ is also estimated, as shown in (15)-(17) at the bottom of the page. The difference between (12)-(14) and (15)-(17) lies in the variable $x_k$: when $x_k = 1$, for lane $j$, it signifies the transformation of the real vehicle $(k, j)$ into a virtual vehicle, while for lane $j'$, it indicates the appearance of a real vehicle $(k, j')$ in this lane, which triggers a change in the state of the virtual vehicle.

$$\eta_2 = \begin{cases} 1, if \ (t_{(p,j')}^d + t_g) mod(t_g + t_r) \neq 0 \\ 0, others \end{cases} \quad (16)$$

The departure time of the vehicles that can receive lane-changing advisory must be within the range of the ATG, as shown in (18) at the bottom of the page.

2) Vehicle lane-changing constrains

The lane-changing maneuver is forbidden for vehicle $(k, j)$ if it stops, as shown in (19).

$$-M \cdot v_{(k,j)} \leq x_k \leq M \cdot v_{(k,j)}, k \in K, j \in J \quad (19)$$

$$t_{(1,j)}^d = \eta_1 \cdot ((1 - x_1) \cdot max\{t_{(p,j)}^d + \tau, t_{(1,j)}^f\} + x_1 \cdot t_{(p,j)}^d) + (1 - \eta_1) \cdot \left((1 - x_1) \cdot t_{(1,j)}^f + x_1 \cdot t_{(p,j)}^d\right), j \in J \quad (12)$$

$$t_{(k,j)}^d = (1 - x_k) \cdot max\{t_{(k-1,j)}^d + \tau, t_{(k,j)}^f, t_G + t_l\} + x_k \cdot t_{(k-1,j)}^d, k \in K \backslash \{k = 1\}, j \in J \quad (14)$$

$$t_{(1,j')}^d = \eta_2 \cdot (x_1 \cdot max\{t_{(p,j')}^d + \tau, t_{(1,j')}^f\} + (1 - x_1) \cdot t_{(p,j')}^d) + (1 - \eta_2) \cdot \left(x_1 \cdot t_{(1,j')}^f + (1 - x_1) \cdot t_{(p,j')}^d\right), j' \in J \quad (15)$$

$$t_{(k,j')}^d = x_k \cdot max\{t_{(k-1,j')}^d + \tau, t_{(k,j')}^f\} + (1 - x_k) \cdot t_{(k-1,j')}^d, k \in K \backslash \{k = 1\}, j' \in J \quad (17)$$

$$min(T_n) \leq t_{(k,j')}^d \leq max(T_n), k \in K, j' \in J, n = 1, 2, \ldots, N \quad (18)$$

$$x_k \cdot \left(d_{safe} + \frac{\left|(v_{(p,j')})^2 - (v_{(k,j)})^2\right|}{2a_L}\right) \leq M \cdot (1 - x_k) + l_{(k,j)} - l_{(p,j')} - l_v, k \in K, j, j' \in J \quad (21)$$

$$x_k \cdot \left(d_{safe} + \frac{\left|(v_{(f,j')})^2 - (v_{(k,j)})^2\right|}{2a_L}\right) \leq M \cdot (1 - x_k) + l_{(f,j')} - l_{(k,j)} - l_v, k \in K, j, j' \in J \quad (22)$$

Vehicles are prohibited from changing lane when they are in close proximity to the stop bar at each intersection along this arterial, which is a common practice in real-world road networks, as shown in (20).

$$l_{(k,j)} - M(1 - x_k) \leq l_c - l_n, k \in K, j \in J \quad (20)$$

3) Lateral safety

(21) and (22) at the bottom of the previous page represent lateral safety constraints, as shown at the bottom of the previous page. When a vehicle decides to change lanes, it must ensure that it maintains a safe distance with both its preceding and following vehicles.

*F. Rolling horizon scheme*

A rolling horizon scheme is proposed for the dynamic implementation of the recognition module and optimization module to adapt to changing traffic conditions, as depicted in Fig. 2. The initialization module and recognition module are performed at each time step. The optimization module is periodically triggered every $h$ seconds. Under the rolling horizon scheme, the control center identifies general traffic in the same lane according to the predetermined plan. It makes appropriate decisions based on real-time vehicle and road information in the next rolling horizon.

IV. NUMERICAL EXAMPLES

*A. Experimental design*

A microsimulation model of a complex intersection is constructed to evaluate the effectiveness of the proposed strategies. Three strategies, namely EBL, BLIDP, and the proposed DSTP, are tested in this simulated environment. The experiments consider scenarios with and without bus stops on the roadway segment. To demonstrate the applicability of the proposed methods, we selected a complex testing scenario, as shown in Fig. 7. Since the proposed method for one direction is independent of the other, this study only considers traffic in one direction. For the driving rules of right-turn vehicles, as depicted in Fig. 7, under the EBL strategy, all right-turn vehicles need to temporarily enter the bus lane from the regular lane before transitioning into the right-turn lane; under the BLIP strategy, only connected right-turn vehicles not within the clearance distance can enter the bus lane, while the rest of the right-turn vehicles follow the same rules as the EBL strategy. It is worth mentioning that in our experiments, the bus lane is located on the far-right side, and left-turning vehicles are not considered within the scope of using the bus lane. Since left-turning vehicles are not considered in allocating road rights and to maximize the flow of through traffic, the experiment includes only through and right-turn vehicles. The surveyed section is approximately 700 m long and extends from three to four lanes, with the transition occurring 100 m from the stop bar. The length of the control area is equal to the length of the surveyed section, i.e., $l_c = 700m$. The length of the no-change zone is $l_n = 30m$. The bus stop is placed approximately 400m from the stop bar for scenarios with the bus stop, i.e., $l_s = 400$. To simulate the randomness of bus stop durations, the dwelling time follows a normal distribution with a mean of 20 seconds and a variance of 10 seconds. The speed limit on the road is $v_d = 13.89$m/s.

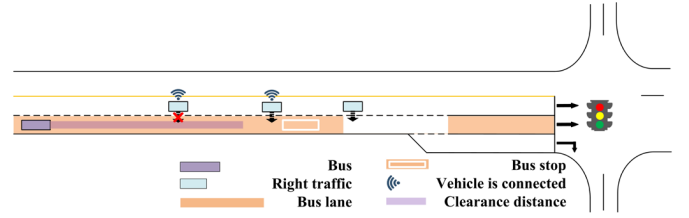

Fig. 7. The geometry layout of simulated roadway

The general traffic, which includes HDVs, CHVs and CAVs, have the same size with a length of 5m and the same performance. The desired acceleration and deceleration values for all general traffic are denoted as 3m/s$^2$ and $-4$m/s$^2$, respectively. The length of the bus is 12m, and the absolute values of the desired acceleration and deceleration for the bus are both denoted as 2m/s$^2$. The CACC model describes the driving behavior of both the bus and CAVs, while the default Krauss model in SUMO describes the driving behavior of HDVs and CHVs. If a CAV follows a HDV or CHV, the car following model of the CAV will switch from CACC to ACC. Due to the unrealistic nature of the default lane-changing model in SUMO, which involves vehicles changing lanes instantaneously, we have replaced it with the sublane model (SL2015) in our experiments. The safe space in our proposed methodology for vehicles to conduct lane-changing is 6m [35]. The maximal physically possible deceleration for general traffic and the bus is $-9$m/s$^2$. The optimization and control interval in the rolling horizon is 5 seconds.

To illustrate the effectiveness of the proposed DSTP strategy, the EBL strategy and the improved BLIDP strategy are applied in this simulation environment as two baselines for comparison. According to previous research, the clear distance for the BLIDP strategy is set to 300m [16]. The arrival of general traffic follows a Poisson distribution. Traffic signals operate under a fixed timing plan, with a cycle length of 100 seconds. Through traffic share the same phase with a green duration of 40 seconds. Signals do not control right-turn vehicles.

We conducted comprehensive experiments on the DSTP strategy, BLIDP strategy, and EBL strategy from four aspects: the connected penetration rates (CPR), traffic demand, bus arrival interval, and right-turn ratio. The CPR in the test scenarios is set at 20%, 40%, 60%, 80%, and 100%. For connected vehicles, the ratio between CAVs and CHVs is 1:1. In addition, the traffic demand in the test scenarios needs to include both under-saturated and over-saturated traffic. The maximum capacity is determined based on the number of vehicles passing through the intersection under the EBL strategy. Since both the DSTP and BLIDP strategies can improve the intersection capacity, six levels of traffic demand are considered (V/C = 0.6, 0.8, 1.0, 1.2, 1.4, 1.6). For scenarios with bus stops, the bus arrival interval follows a normal distribution with a mean of 30, 60, 90, 120, 180, and 240 seconds and a standard deviation of 10 seconds. For scenarios without bus stops follows a normal distribution with a mean of 20, 30, 60, 90, 120, 180, and 240 seconds and a standard deviation of 10 seconds. Different proportions of right-turn vehicles (0.1, 0.3, 0.5, 0.7, 0.9) are also tested in the simulation.

The EBL, BLIDP, and proposed DSTP strategy are implemented in Python. The optimization model is solved using

Gurobi 10.0. All simulation work is conducted using SUMO [36] on a desktop computer with a 4.9-GHz Intel Core processor and 32 GB of RAM. The default lane-changing models in SUMO simulate the driving behavior of general traffic. It is worth noting that, to fully replicate real-world traffic conditions, although lane-changing advisories are sent to vehicles in advance through the control center, the vehicles are not forcibly required to change lanes. Instead, they autonomously determine whether they can change lanes based on the lane-changing models. Five random seeds are used in the simulation and each simulation runs for 1800 s with a full warm-up period. The resolution of the simulator is 0.5 s.

*B. Results and discussions*

In each section, we discussed the performance of through general traffic, right-turn general traffic, and buses in terms of delay, fuel consumption, and $CO_2$ emissions. The Handbook Emission Factors for Road Transport model (HBEFA) [37], provided in SUMO, is used to estimate fuel consumption and $CO_2$ emissions.

*1) Impact of the connected penetration rates*

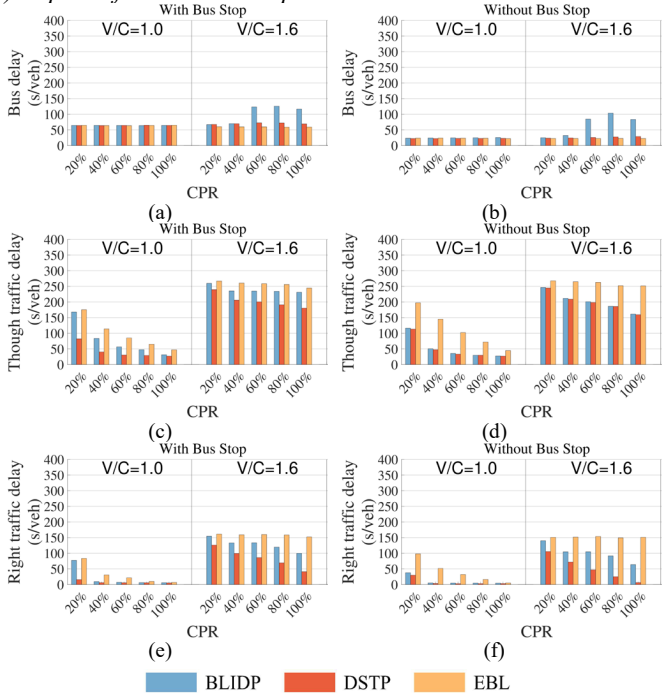

Fig. 8. Delay of bus and general traffic under varying CPR: (a) bus delay in the bus stop scenario; (b) bus delay without the bus stop scenario; (c) through traffic delay in the bus stop scenario; (d) through traffic delay without the bus stop scenario; (e) right-turn traffic delay in the bus stop scenario; (f) right-turn traffic delay without the bus stop scenario

In this section, we tested the performance of the DSTP, BLIDP, and EBL strategies under different CPR. In this case, bus arrival interval is 60s and the right-turn ratio is 0.3. To demonstrate the performance of the compared strategies more precisely under medium and high traffic demand conditions, we conducted tests with V/C values of 1 and 1.6, respectively. Fig. 8 illustrates the impact of the three strategies on bus and general traffic delays as CPR increases, both with and without bus stops.

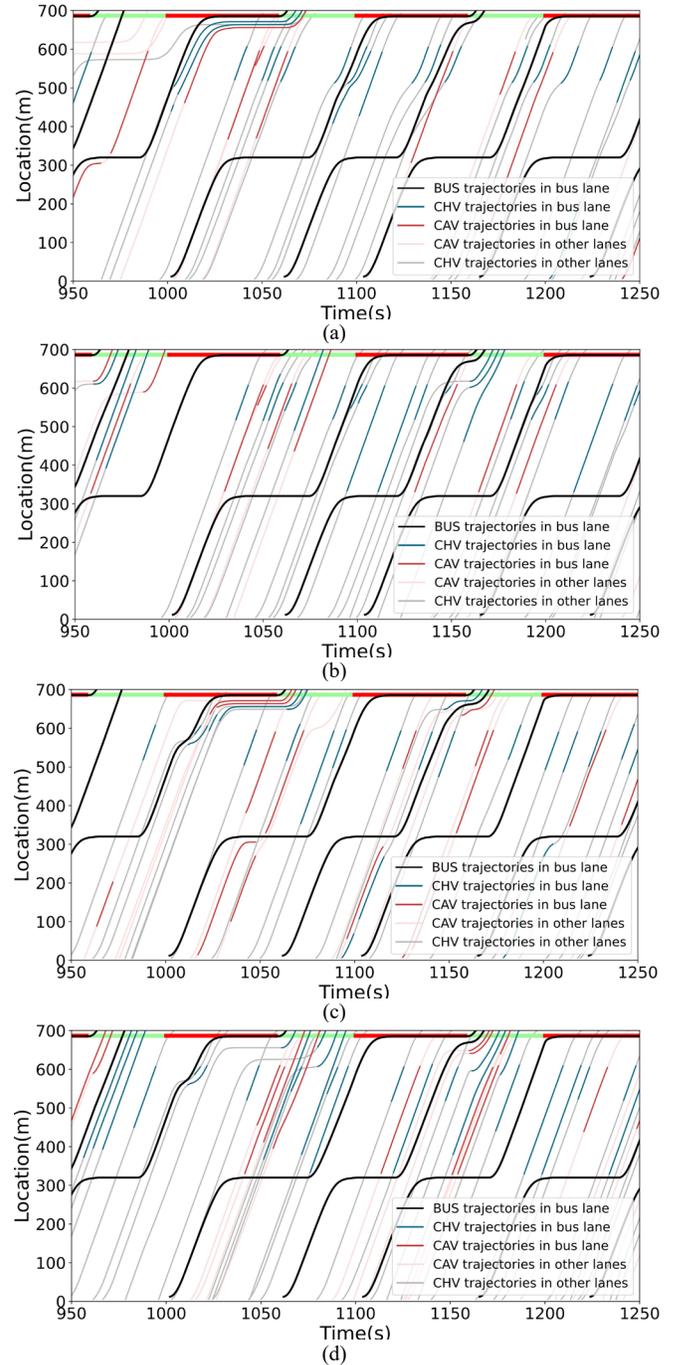

Fig. 9. Spatial-temporal trajectories of vehicles in bus lane (V/C = 1): (a) BLIDP strategy, CPR = 40%; (b) DSTP strategy, CPR = 40%; (c) BLIDP strategy, CPR =80%; (d) DSTP strategy, CPR =80%.

It can be observed that at V/C=1, the BLIDP strategy and the DSTP strategy ensure priority for buses at any CPR. In this case, the advantage of the DSTP strategy mainly lies in reducing delays for through and right-turn traffic in scenarios with bus stop, as shown in Fig. 8(c) and (e). At low CPR, the DSTP strategy significantly outperforms the BLIDP strategy in reducing delays. As CPR increases, although the effectiveness of the BLIDP strategy gradually approaches that of the DSTP strategy, the DSTP strategy remains superior. This is because when buses stop, the rules of the BLIDP strategy prevent general traffic within a certain clearance distance from entering the bus lane, especially when low CPR situation and clearance distance overlap, resulting in fewer general traffic vehicles able to enter the bus lane and leading to resource wastage, as shown in the trajectory plots in Fig. 9 (a) and (c). Conversely, under the DSTP strategy, buses stopping actually encourage the utilization of remaining resources in the bus lane by general traffic, as illustrated in Fig. 9 (b) and (d). Therefore, compared to the other two strategies, the DSTP strategy is more applicable in ensuring bus priority and reducing delays for general traffic, even at low CPR.

Furthermore, at a higher traffic demand level (V/C=1.6), the advantage of the DSTP strategy lies mainly in ensuring stable bus delay under general traffic borrowing bus lane conditions and reducing delays for right-turn traffic. From Fig. 8(a) and (b), compared to the BLIDP strategy, the increase in bus delay under the DSTP strategy is not significant as CPR increases. Moreover, from Fig. 8(e) and (f), it is evident that the delay for right-turn traffic under the DSTP strategy is significantly lower than that under the BLIDP strategy.

Fig. 10 further explains the reasons for the above phenomena. For the BLIDP strategy, at low CPR, as shown in Fig. 10 (a), some general traffic chooses to enter the bus lane as soon as they enter the road. However, when the bus arrives, the general traffic density on the general lanes increases, leading to more general traffic being unable to find space for lane changing to exit the bus lane. These vehicles have to continue driving in the bus lane, causing queues even at bus stops. As CPR increases, as shown in Fig. 10 (c), more general traffic can enter the bus lane, leading to further congestion of general traffic that cannot be cleared in time, resulting in queues of vehicles before the stop bar and causing bus delays. At the same time, some right-turn traffic also needs to queue with the through traffic ahead, resulting in significant delays for right-turn traffic. In contrast, as shown in Fig. 10 (b) and (d), under the DSTP strategy, bus delays are not significantly affected by the high CPR. This is because the DSTP strategy optimally combines green time with available lane resources to organize general traffic in the bus lane more efficiently. In other words, the DSTP strategy selects appropriate vehicles to enter the bus lane through optimization methods, avoiding severe queues of general traffic on the bus lane.

We further depicted the fuel economy and CO2 emissions of vehicles under bus arrival interval is 60s and V/C is 1.6 conditions. Fig. 11 illustrates the fuel economy of vehicles under three strategies at different traffic demand levels. It can be observed that with CPR increases, the fuel economy of buses under the BLIDP strategy decreases to a certain extent, while the fuel economy of buses under the DSTP strategy remains similar to that under the EBL strategy, showing no significant changes. Furthermore, in the presence of bus stops, as CPR increases, both the BLIDP and EBL strategies lead to a decrease in the fuel economy of through traffic, while the fuel economy of through traffic under the DSTP strategy does not change significantly. Lastly, a clear advantage of the DSTP strategy is that with the increase in CPR, the fuel economy of right-turn traffic gradually increases and is significantly higher than that under the BLIDP strategy.

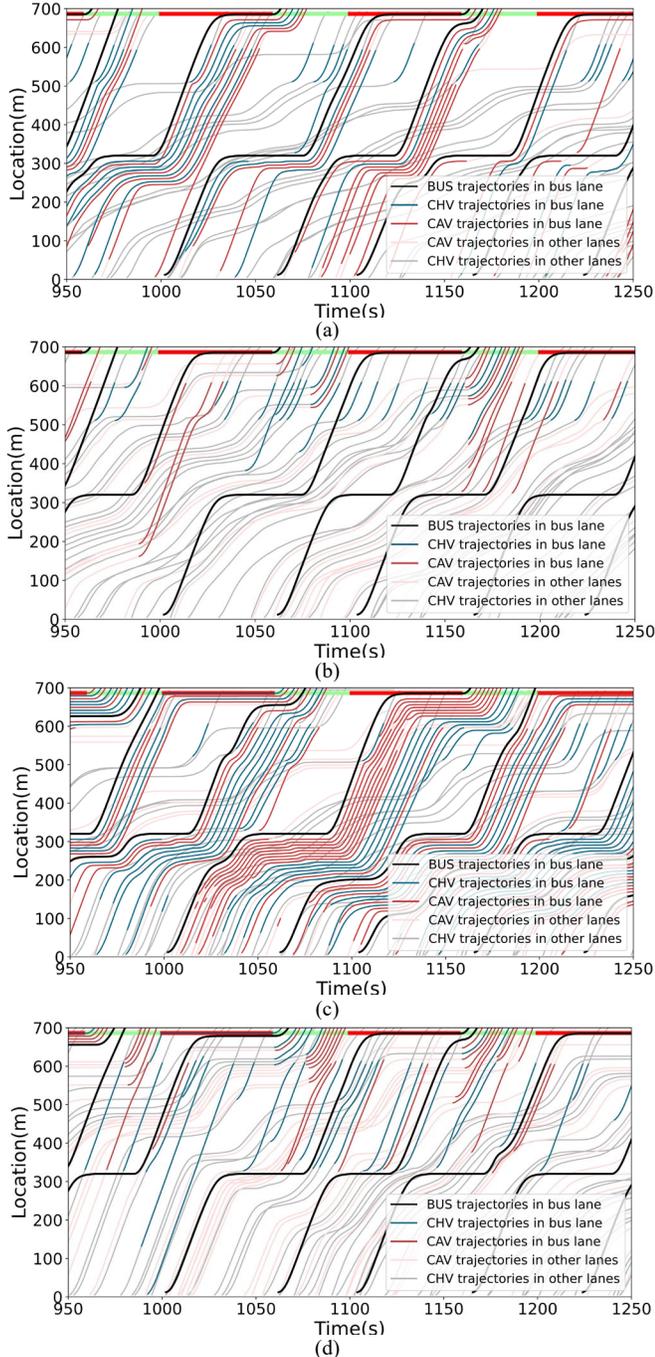

Fig. 10. Spatial-temporal trajectories of vehicles in bus lane (V/C = 1.6): (a) BLIDP strategy, CPR = 40%;(b) DSTP strategy, CPR = 40%;(c) BLIDP strategy, CPR =80%;(d) DSTP strategy, CPR =80%.

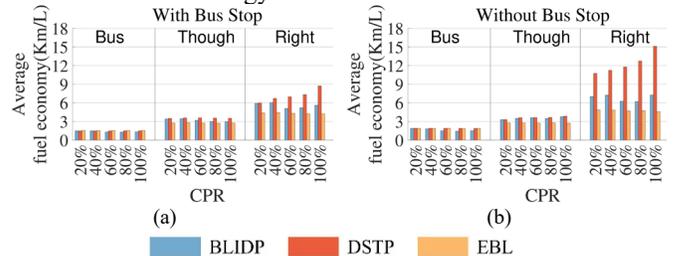

Fig. 11. Impacts on fuel economy under varying CPR: (a) in the bus stop scenario; (b) without the bus stop scenario

Fig. 12 presents a stacked bar chart of the average CO2 emissions of buses, through traffic, and right-turn traffic under the three strategies at different CPR. It can be seen that the DSTP strategy is better at minimizing the CO2 emissions of all vehicles on the road, especially when the CPR exceeds 60%.

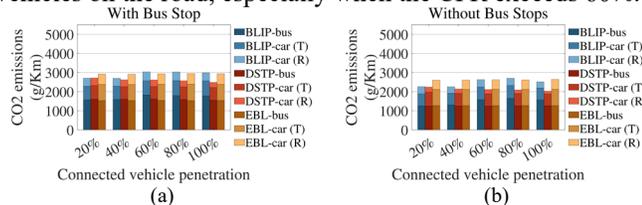

Fig. 12. Impacts on CO2 emissions under varying CPR: (a) fuel economy in the bus stop scenario; (b) fuel economy without the bus stop scenario

*2) Impact of traffic volume*

In this section, we tested the impact of different traffic volumes on the three strategies. Considering that both the BLIDP and DSTP strategies performed best at 100% CPR, to further highlight the differences between DSTP and BLIDP, we set CPR to 100% in this experiment, with the penetration rate of CAVs set to 50% to emphasize the complexity of mixed traffic flow and the necessity of Right-of-Way assignment optimization. Moreover, to comprehensively demonstrate the influence of traffic volume, three levels of bus arrival intervals (30s, 60s, 90s) for each scenario were also presented.

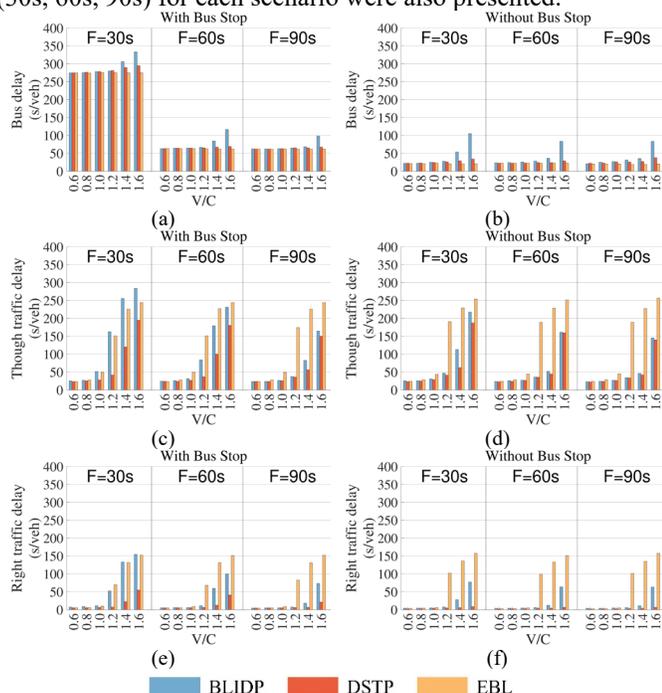

Fig. 13. Impacts on delay of bus and general traffic under varying traffic demand levels: (a) bus delay in the bus stop scenario; (b) bus delay without the bus stop scenario; (c) through traffic delay in the bus stop scenario; (d) through traffic delay in the bus stop scenario; (e) right-turn traffic delay in the bus stop scenario; (f) right-turn traffic delay without the bus stop scenario

Fig. 13 (a), (c), and (e) depict the average delays of buses, through traffic, and right-turn traffic with a bus stop, whereas Fig. 13 (b), (d), and (f) illustrate the corresponding delays without a bus stop, both under varying traffic demand levels. It can be seen that the BLIDP strategy is significantly influenced by the traffic demand level. The DSTP strategy demonstrates superior effectiveness in ensuring bus priority and reducing general traffic delays under high-demand conditions. Firstly, for buses, both from Fig. 13 (a) and (b), it is evident that as V/C increases, the rise in bus delays under the DSTP strategy is significantly lower compared to the BLIDP strategy, and is relatively close to the delays under the EBL strategy. Taking the scenario of V/C=1.6 and bus arrival interval is 30s without bus stops as an example, under the BLIDP strategy, bus delays are 84 seconds higher than under the EBL strategy, while under the DSTP strategy, bus delays are only 13 seconds higher than under the EBL strategy, indicating that the DSTP strategy is more effective in ensuring bus priority under high V/C conditions. Secondly, for general traffic, whether under bus arrival interval is 30s, 60s, or 90s, with or without bus stops, the effectiveness of the DSTP strategy in reducing delays increases compared to the BLIDP strategy as V/C rises. Furthermore, from Fig. 13 (e) and (f), it can be observed that the right-turn traffic delays under the DSTP strategy are significantly lower than under the other two strategies. Taking the scenario without bus stops at bus arrival interval is 90s and V/C=1.6 as an example, although the DSTP and BLIDP strategies have comparable performance in reducing through traffic delays, the right-turn traffic delays under the DSTP strategy are 60 seconds lower than under the BLIDP strategy.

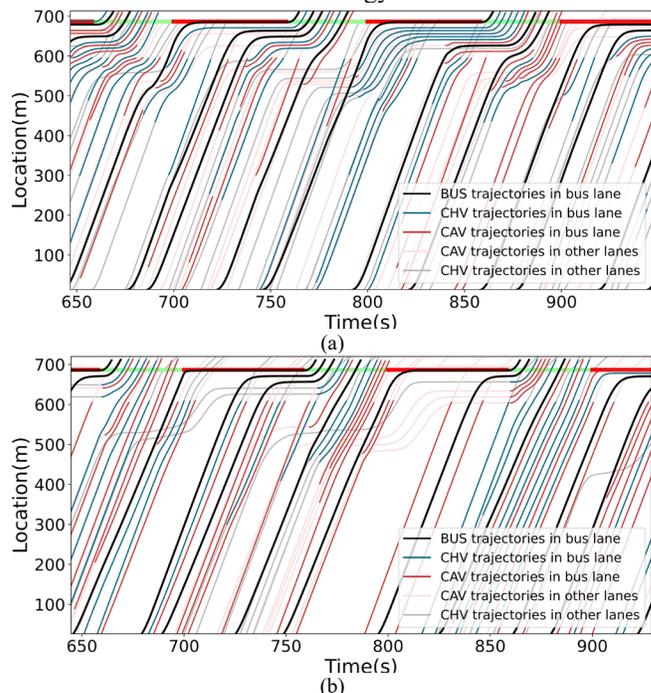

Fig. 14. Spatial-temporal trajectories of vehicles in bus lane without the bus stop scenario: (a) BLIDP strategy, V/C=1.4; (b) DSTP strategy, V/C=1.4

Fig. 14 provides a more visual representation of the reasons behind the aforementioned results. Fig. 14 (a) shows the vehicle trajectories in the bus lane under the BLIDP strategy when V/C=1.4 and bus arrival interval is 30s. It can be seen that the operations of the second, sixth, and eighth buses are significantly disrupted because vehicles that cannot exit the bus lane quickly form queues before the stop bar. Although the third and seventh buses successfully clear the general traffic ahead, a significant portion of green time is not utilized. More importantly, they would arrive at the stop bar during red time even if they cruise without disruptions, and clearing the vehicles ahead would only waste spatial-temporal resources in the bus lane. Additionally, the BLIDP strategy may lead to long queues of general traffic in the bus lane, hindering the smooth entry of right-turn traffic into the right-turn lane. For instance, during the first green time in Fig. 14 (a), right-turn vehicles are unable to enter the right-turn lane promptly due to the presence of queueing vehicles in the bus lane. In the same scenario, as shown in Fig. 14 (b), the DSTP strategy provides a more optimal bus lane sharing scheme. It can be observed that although general traffic travels in the bus lane, they have minimal impact on the smooth flow of buses. However, due to the stochastic nature of vehicle movements, the prediction of travel times is not entirely accurate, and some general traffic may slightly affect the smooth cruising of buses, which is also the reason for the increase in bus delays shown in Fig. 13(a) and (b) under the DSTP strategy as V/C increases. For example, after the end of the third green time in Fig. 14 (b), two vehicles remain before the stop bar in the bus lane, causing slight fluctuations in the trajectory of the following bus. Moreover, right-turn traffic operates more smoothly under the DSTP strategy because this strategy does not lead to significant queues in the bus lane. Consequently, the operation of right-turn traffic is nearly unimpeded and can smoothly enter the right-turn lane from the bus lane.

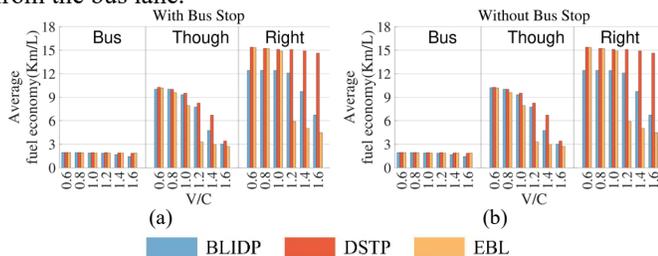

Fig. 15. Impacts on fuel economy under varying traffic demand levels: (a) in the bus stop scenario; (b) without the bus stop scenario

We further examined the performance of three strategies in terms of fuel economy and CO2 emissions under dense bus arrival scenario as traffic demand levels increase. The bus arrival interval and right-turn ratio are set as 30 and 0.3, respectively. Fig. 15 shows the average fuel economy of buses, through traffic, and right-turn traffic under various traffic demand levels. Observing Fig. 15, it is evident that the fuel economy of buses remains relatively consistent and unaffected by the traffic demand level under both the EBL and DSTP strategies. However, the BLIDP strategy shows a gradual decrease in fuel economy of buses as the traffic demand level increases, although it still maintains a similar level to the other two strategies when V/C<1. Furthermore, while the fuel economy of through and right-turn traffic decreases with increasing traffic demand levels under all three strategies, the advantages of the DSTP strategy become evident when V/C≥1, especially in scenarios with bus stop. Another noteworthy result is the significant impact of the presence or absence of bus stop on the effectiveness of the BLIDP strategy. In scenarios without bus stops, the fuel economy of the BLIDP strategy is higher than that of the EBL strategy, whereas in scenarios with bus stops, the fuel economy of the BLIDP strategy is even the lowest among the three strategies. Therefore, in terms of average fuel economy, the DSTP strategy provides a more substantial fuel savings advantage than the other two strategies. It not only maintains the fuel economy of buses at a similar level to the EBL strategy but also ensures optimal fuel economy for through traffic compared to the other two strategies. Additionally, the DSTP strategy maintains high fuel efficiency for right-turn traffic, even under high traffic demand levels.

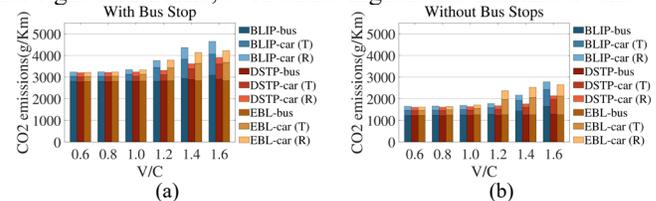

Fig. 16. Impacts on CO2 emissions under varying traffic demand levels: (a) in the bus stop scenario; (b) without the bus stop scenario

Fig. 16 illustrates the average CO2 emissions of buses, through traffic, and right-turn traffic under different traffic demand levels. It can be observed that although the EBL strategy reduces CO2 emissions for buses, it increases CO2 emissions for general traffic. While the BLIDP strategy controls CO2 emissions for general traffic to some extent compared to the EBL strategy, as the traffic demand level increases, CO2 emissions for general traffic gradually become uncontrolled, and CO2 emissions for buses also increase. The DSTP strategy proves to be a suitable choice as it effectively controls CO2 emissions for buses under any demand level, and it performs better in controlling CO2 emissions for general traffic compared to the BLIDP strategy, particularly in situations where traffic demand levels are high.

*3) Impact of bus arrival interval*

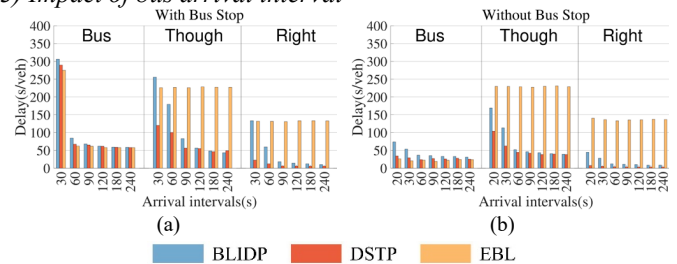

Fig. 17. Delay of bus and general traffic under different bus arrival intervals: (a) in the bus stop scenario; (b) without the bus stop scenario

Fig. 17 illustrates the influence of varying bus arrival intervals on delays when implementing different strategies under fully connected environment. In this case, V/C=1.4, right-ratio = 0.3. It is evident that under the EBL strategy, the delay for buses and general traffic is almost unaffected by the bus arrival interval. It is worth noting that in Fig. 17 (a), the significant increase in bus delay when bus arrival interval = 30s is attributed to the high bus arrival interval and the necessity for buses to stop sequentially at the bus stop, leading to queues forming before the stop bar and bus stops. Furthermore, both the BLIDP and DSTP strategies have a positive effect on reducing general traffic delays. Comparatively, as bus arrivals become more frequent, the DSTP strategy exhibits a more pronounced reduction in general traffic delays and causes lower delays for buses compared to the BLIDP strategy. This is because the BLIDP strategy's rules prohibit general traffic from entering within a certain distance ahead of the bus when buses are stopping at this bus continuously, resulting in underutilization of the vacant time and space. In contrast, the DSTP strategy disregards the clearance distance constraint and dynamically identifies the empty time and space in the bus lane. Therefore, the DSTP strategy is more effective than the BLIDP strategy in reducing general traffic delays while prioritizing buses, especially in dense bus arrivals. As bus arrivals become increasingly sparse, the advantages of both strategies in reducing general traffic delays converge, and their effect on bus delays gradually aligns with that of the EBL strategy. This is due to the diminishing number of buses arriving, causing the bus lane under the BLIDP and DSTP strategies to gradually resemble a general traffic lane.

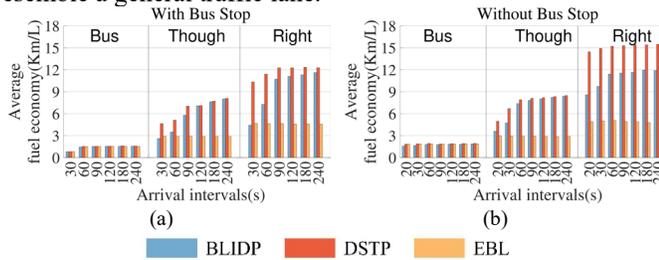

Fig. 18. Fuel economy of bus and general traffic under different bus arrival intervals: (a) in the bus stop scenario; (b) without the bus stop scenario

Fig. 18 illustrates the influence of varying bus arrival intervals on fuel economy when implementing different strategies. It can be observed that as bus arrival intervals increase, although the fuel economy of general traffic under both the BLIDP and DSTP strategies gradually improves, the DSTP strategy exhibits higher fuel economy, especially in situations with dense bus arrivals. Additionally, it is worth noting that in scenarios without bus stops, the DSTP strategy ensures higher fuel economy for right-turn traffic compared to the BLIDP strategy.

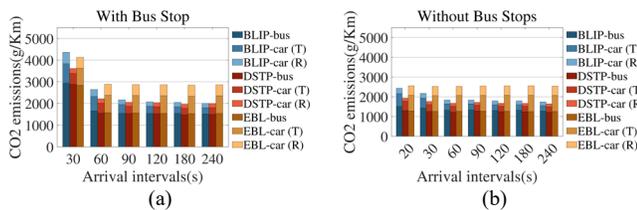

Fig. 19. $CO_2$ emissions of bus and general traffic under different bus arrival intervals: (a) in the bus stop scenario; (b) without the bus stop scenario

Fig. 19 illustrates the influence of varying bus arrival intervals on $CO_2$ emissions when implementing different strategies. Comparing the DSTP strategy with the EBL strategy, it can be observed that under the DSTP strategy, $CO_2$ emissions for buses are similar to those under the EBL strategy, while $CO_2$ emissions for general traffic are significantly lower than those under the EBL strategy. Furthermore, comparing the DSTP strategy with the BLIDP strategy, it is evident that under dense bus arrivals, both $CO_2$ emissions for buses and general traffic are lower under the DSTP strategy compared to the BLIDP strategy. As the bus arrival interval increases, the $CO_2$ emissions for vehicles under both strategies tend to be the same.

*4) Impact of right-turn ratios*

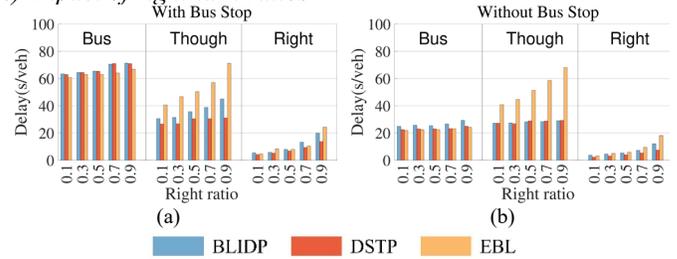

Fig. 20. Delay of bus and general traffic under different right ratios: (a) in the bus stop scenario; (b) without the bus stop scenario

Fig. 20 illustrates the influence of varying right-turn ratios on delays when implementing different strategies under fully connected environment, with V/C at 1.0 and a bus arrival interval of 60s. As the demand of right-turn traffic increases, the DSTP strategy demonstrates superiority over the BLIDP strategy in reducing delays for right-turn traffic. This is attributed to the DSTP strategy prioritizing right-turn traffic by utilizing the bus lane, whereas the BLIDP strategy does not explicitly address the priority of right-turn traffic. Consequently, under the BLIDP strategy, through traffic already in the bus lane may impede right-turn traffic in the adjacent lane, making it difficult for right-turn traffic to promptly find adequate space to transition into the bus lane. Furthermore, compared to the BLIDP strategy, the advantage of the DSTP strategy in scenario with bus stop lies in reducing delays for through traffic, while in scenarios without bus stops, the DSTP strategy excels in ensuring that bus delays do not significantly increase. Overall, the DSTP strategy leads to more consistent delays for buses and general traffic, showing reduced sensitivity to the rise in right-turn traffic.

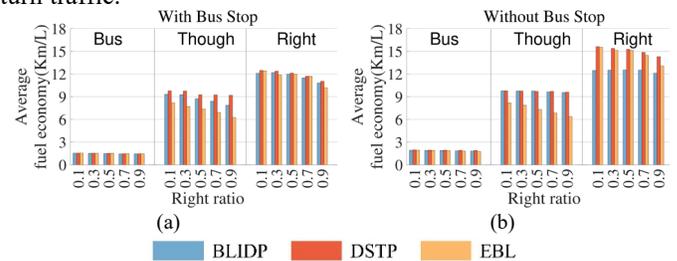

Fig. 21. Fuel economy of bus and general traffic under different right ratios: (a) in the bus stop scenario; (b) without the bus stop scenario

Fig. 21 illustrates the influence of varying right ratios on fuel economy when implementing different strategies. For buses, their fuel economy is less influenced by the right-turn ratios, regardless of the strategy employed. For through traffic, the EBL strategy is more sensitive to the right-turn ratios, resulting in lower fuel economy compared to the other two strategies. While both the BLIDP and DSTP strategies demonstrate less sensitivity to the right-turn ratios regarding fuel economy, the DSTP strategy consistently surpasses the BLIDP strategy in enhancing fuel economy, especially in scenarios with bus stop. For right-turn traffic, although the fuel economy decreases as the right-turn ratio increases under all three strategies, the DSTP strategy still achieves higher fuel economy than the other two strategies, especially in the absence of a bus stop scenario.

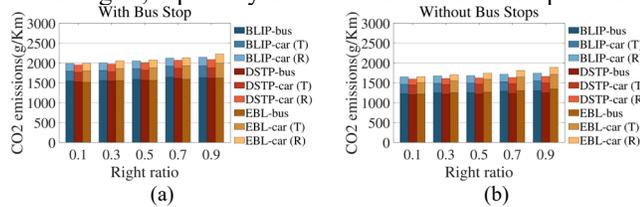

Fig. 22. CO2 emissions of bus and general traffic under different right ratios: (a) in the bus stop scenario; (b) without the bus stop scenario

Fig. 22 illustrates the influence of varying right-turn ratios on CO2 emissions when implementing different strategies. It can be observed that as the right-turn ratio increases, the CO2 emissions for all three strategies show an upward trend. However, the EBL strategy exhibits the largest increase, followed by the BLIDP strategy, while the DSTP strategy ensures the lowest rise in CO2 emissions.

## V. Conclusion

This study proposes an active control strategy for bus lanes termed Dynamic Spatial-Temporal Priority (DSTP) under partially connected vehicle environment, with the objective of optimizing the utilization of the residual spatial-temporal resources within the bus lane for both CHVs and CAVs, while concurrently ensuring priority for bus operations. The DSTP strategy involves the estimation of vehicle departure times from the stop bar, considering individual vehicles based on their longitudinal locations. Subsequently, the control center identifies remaining resources within the bus lane based on the departure times, speeds, and locations of vehicles. Using the proposed Right-of-Way assignment optimization model, eligible private cars are selected from general lanes to enter the bus lane to alleviate congestion on general lanes. The implementation of this strategy incorporates a rolling horizon scheme to dynamically adapt to time-varying traffic conditions. This allows for continuous adjustment of the strategy based on the real-time traffic condition.

Simulation experiments validate the effectiveness of the proposed DSTP strategy and compares it with the EBL and BLIDP strategies across various metrics, including delay, fuel economy, and CO2 emissions. The results indicate that the DSTP strategy integrates the advantages of EBL and BLIDP, maintaining bus priority akin to EBL while enhancing traffic efficiency and reducing fuel consumption, and CO2 emissions, especially during high traffic demand and concentrated bus arrivals. Although this paper is centered on bus priority, the proposed strategy applies to other priority vehicles. For example, it can be extended to scenarios where CAVs have priority on dedicated CAV lanes while allowing other vehicles to enter. In conclusion, this strategy can improve the traffic efficiency of roads in mixed traffic environments with bus lanes and reduce environmental pollution.

The study's assumption of fixed signal timing plans at the intersection, constraining potential operational performance improvements. Therefore, future researches need to integrate signal timing optimization with the DSTP strategy in a mixed traffic environment. Additionally, the study assumes that all vehicles receiving lane change advisories will attempt to change lanes unless the current driving conditions present challenges for them to safely making lane-changes. However, in reality, not all drivers who receive lane change advisories will attempt to change lanes. Therefore, future research needs to consider driver compliance with advisories in various real-world traffic scenarios.

## References

[1] L.T. Truong, G. Currie, M. Wallace and C. De Gruyter, "Analytical Approach to Estimate Delay Reduction Associated with Bus Priority Measures," *Ieee Intell. Transp. Syst. Mag.*, vol. 9, no. 4, pp. 91-101, 2017.

[2] H. Al-Deek, A. Sandt, A. Alomari and O. Hussain, "A technical note on evaluating the effectiveness of bus rapid transit with transit signal priority," *J. Intell. Transport. Syst.*, vol. 21, no. 3, pp. 227-238, 2017.

[3] A. Alam, E. Diab, A.M. El-Geneidy and M. Hatzopoulou, "A simulation of transit bus emissions along an urban corridor: Evaluating changes under various service improvement strategies," *Transportation Research Part D: Transport and Environment*, vol. 31, pp. 189-198, 2014, doi: https://doi.org/10.1016/j.trd.2014.06.010.

[4] D. Tsitsokas A. Kouvelas and N. Geroliminis, "Modeling and optimization of dedicated bus lanes space allocation in large networks with dynamic congestion," *Transportation Research Part C: Emerging Technologies*, vol. 127, pp. 103082, 2021.

[5] X. Bai, Z. Zhou, K. Chin and B. Huang, "Evaluating lane reservation problems by carbon emission approach," *Transportation Research Part D: Transport and Environment*, vol. 53, pp. 178-192, 2017, doi: https://doi.org/10.1016/j.trd.2017.04.002.

[6] J. Viegas and B. Lu, "Widening the scope for bus priority with intermittent bus lanes," *Transp. Plan. Technol.*, vol. 24, no. 2, pp. 87-110, 2001.

[7] J. Viegas and B. Lu, "The Intermittent Bus Lane signals setting within an area," *Transportation Research Part C: Emerging Technologies*, vol. 12, no. 6, pp. 453-469, 2004.

[8] M. Eichler and C.F. Daganzo, "Bus lanes with intermittent priority: Strategy formulae and an evaluation," *Transportation Research Part B: Methodological*, vol. 40, no. 9, pp. 731-744, 2006.

[9] N. Chiabaut X. Xie and L. Leclercq, "Road Capacity and Travel Times with Bus Lanes and Intermittent Priority Activation," *Transportation Research Record: Journal of the Transportation Research Board*, vol. 2315, no. 1, pp. 182-190, 2012.

[10] S.I. Guler and M.J. Cassidy, "Strategies for sharing bottleneck capacity among buses and cars," *Transportation Research Part B: Methodological*, vol. 46, no. 10, pp. 1334-1345, 2012.

[11] G. Carey T. Bauer and K. Giese, Bus lane with intermittent priority (BLIMP) concept simulation analysis final report : November 2009., *Book* Bus lane with intermittent priority (BLIMP) concept simulation analysis final report : November 2009., Series Bus lane with intermittent priority (BLIMP) concept simulation analysis final report :

November 2009.,ed., Editor ed., 2009, pp.

[12] F. Qiu, W. Li, J. Zhang, X. Zhang and Q. Xie, "Exploring suitable traffic conditions for intermittent bus lanes," *J. Adv. Transp.*, vol. 49, no. 3, pp. 309-325, 2015.

[13] H.B. Zhu, "Numerical study of urban traffic flow with dedicated bus lane and intermittent bus lane," *Physica a: Statistical Mechanics and its Applications*, vol. 389, no. 16, pp. 3134-3139, 2010.

[14] C.X. Ma and X.C. Xu, "Providing Spatial-Temporal Priority Control Strategy for BRT Lanes: A Simulation Approach," *J. Transp. Eng. Pt a-Syst.*, vol. 146, no. 7, 2020.

[15] K. Othman A. Shalaby and B. Abdulhai, "Dynamic Bus Lanes Versus Exclusive Bus Lanes: Comprehensive Comparative Analysis of Urban Corridor Performance," *Transportation Research Record: Journal of the Transportation Research Board*, pp. 862748365, 2022.

[16] W. Wu, L. Head, S. Yan and W. Ma, "Development and evaluation of bus lanes with intermittent and dynamic priority in connected vehicle environment," *J. Intell. Transport. Syst.*, vol. 22, no. 4, pp. 301-310, 2018.

[17] Y. Feng, K.L. Head, S. Khoshmagham and M. Zamanipour, "A real-time adaptive signal control in a connected vehicle environment," *Transportation Research Part C: Emerging Technologies*, vol. 55, pp. 460-473, 2015.

[18] Z. Yang, Y. Feng and H.X. Liu, "A cooperative driving framework for urban arterials in mixed traffic conditions," *Transportation Research Part C: Emerging Technologies*, vol. 124, pp. 102918, 2021.

[19] M.W. Levin and A. Khani, "Dynamic transit lanes for connected and autonomous vehicles," *Public Transport*, vol. 10, no. 3, pp. 399-426, 2018.

[20] D. Ou, R. Liu, I. Rasheed, L. Shi and H. Li, "Operation performance of tram lanes with intermittent priority with the coexistence of regular and automatic vehicles," *J. Intell. Transport. Syst.*, vol. 26, no. 4, pp. 486-497, 2022.

[21] Y. Luo, J. Chen, S. Zhu and Y. Yang, "Developing the Dynamic Bus Lane Using a Moving Block Concept," *Transportation Research Record: Journal of the Transportation Research Board*, pp. 862748314, 2022.

[22] A. Kampouri and I. Politis, Optimization of a Bus Lane with Intermittent Priority dynamically activated by the road traffic, *Proc. In Proceedings of the 23rd International Transport and Air Pollution Conference*, 2020, pp. 314.

[23] X. Chen, X. Lin, F. He and M. Li, "Modeling and control of automated vehicle access on dedicated bus rapid transit lanes," *Transportation Research Part C: Emerging Technologies*, vol. 120, pp. 102795, 2020.

[24] P. Anderson and N. Geroliminis, "Dynamic lane restrictions on congested arterials," *Transportation Research Part a: Policy and Practice*, vol. 135, pp. 224-243, 2020.

[25] Y. Shao, J. Sun, Y. Kan and Y. Tian, "Operation of dedicated lanes with intermittent priority on highways: conceptual development and simulation validation," *J. Intell. Transport. Syst.*, vol. ahead-of-print, no. ahead-of-print, pp. 1-15, 2022.

[26] X. Shan, C. Wan, P. Hao, G. Wu and M.J. Barth, "Developing A Novel Dynamic Bus Lane Control Strategy With Eco-Driving Under Partially Connected Vehicle Environment," *Ieee Trans. Intell. Transp. Syst.*, pp. 1-16, 2024.

[27] G. Zheng and M. Zhang, Intermittent Bus Lane Control Method for Preventing Blocking of Right-Turning Vehicles, IEEE, 2021, pp. 4022-4027.

[28] V. Milanes, S.E. Shladover, J. Spring, C. Nowakowski, H. Kawazoe and M. Nakamura, "Cooperative Adaptive Cruise Control in Real Traffic Situations," *Ieee Trans. Intell. Transp. Syst.*, vol. 15, no. 1, pp. 296-305, 2014.

[29] V. Milanés and S.E. Shladover, "Modeling cooperative and autonomous adaptive cruise control dynamic responses using experimental data," *Transportation Research Part C: Emerging Technologies*, vol. 48, pp. 285-300, 2014.

[30] C. Zhang, N.R. Sabar, E. Chung, A. Bhaskar and X. Guo, "Optimisation of lane-changing advisory at the motorway lane drop bottleneck," *Transportation Research Part C: Emerging Technologies*, vol. 106, pp. 303-316, 2019.

[31] S. Ilgin Guler, M. Menendez and L. Meier, "Using connected vehicle technology to improve the efficiency of intersections," *Transportation Research Part C: Emerging Technologies*, vol. 46, pp. 121-131, 2014.

[32] X.J. Liang, S.I. Guler and V.V. Gayah, "An equitable traffic signal control scheme at isolated signalized intersections using Connected Vehicle technology," *Transportation Research Part C: Emerging Technologies*, vol. 110, pp. 81-97, 2020.

[33] C. Ma, C. Yu and X. Yang, "Trajectory planning for connected and automated vehicles at isolated signalized intersections under mixed traffic environment," *Transportation Research Part C: Emerging Technologies*, vol. 130, pp. 103309, 2021.

[34] Z. Zheng, "Recent developments and research needs in modeling lane changing," *Transportation Research Part B: Methodological*, vol. 60, pp. 16-32, 2014, doi: 10.1016/j.trb.2013.11.009.

[35] Q. Wang, Y. Gong and X. Yang, "Connected automated vehicle trajectory optimization along signalized arterial: A decentralized approach under mixed traffic environment," *Transportation Research Part C: Emerging Technologies*, vol. 145, pp. 103918, 2022.

[36] D. Krajzewicz, J. Erdmann, M. Behrisch and L. Bieker, "Recent development and applications of SUMO-Simulation of Urban MObility," *International Journal On Advances in Systems and Measurements*, vol. 5, no. 3&4, 2012.

[37] S. Hausberger, "Emission Factors from the Model PHEM for the HBEFA Version 3,", 2009.

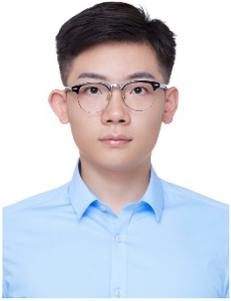

Haoran Li (Graduate Student Member, IEEE) received the B.S. degree in transportation engineering from Shijiazhuang Tiedao University, Shijiazhuang, China. He is currently pursuing the Ph.D. degree in School of Traffic and Transportation, Beijing Jiaotong University, Beijing, China. His research interests include bus priority strategy, intelligent transportation systems, connected automated vehicles, and transportation simulation.

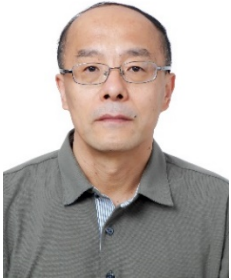

Zhenzhou Yuan received Ph.D. degree in transportation planning and management from Beijing Jiaotong University, Beijing, China, in 2000. He is currently a Professor with the School of Traffic and Transportation, Beijing Jiaotong University, Beijing, China. His research interests include transportation systems optimization, traffic safety analysis and intelligent transportation system (ITS) technology development.

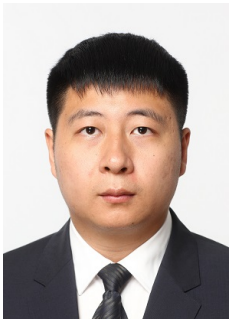

Rui Yue is an assistant professor at the Beijing Jiaotong University, Beijing, China. He completed his Bachelor's degree in Civil Engineering in 2015 at Beijing Jiaotong University, followed by a Master's degree in Transportation Engineering in 2017 at the University of Nevada, Reno. In 2020, he earned his Ph.D. in Civil and Environmental Engineering from the University of Nevada, Reno.

His primary research areas encompass traffic operations and signal control, capacity modeling and driver behavior analysis based on trajectories performance assessment for alternative intersections and interchanges, and LiDAR applications. Dr. Yue serves as a young academic editor for the Journal of "Digital Transportation and Safety" and "Municipal Engineering Technology".

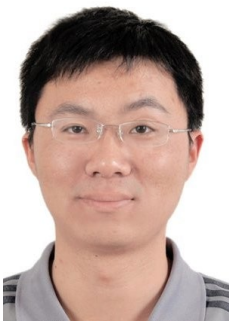

Guangchuan Yang is a Research Associate at the Institute for Transportation Research and Education at North Carolina State University. He completed his Bachelor's degree in Traffic Engineering in 2010 at Dalian Jiaotong University, followed by a Master's degree in Transportation Planning and Engineering in 2012 at the University of Southampton. In 2017, he earned his Ph.D. in Civil and Environmental Engineering from the University of Nevada, Reno.

His primary research areas encompass traffic operations and safety modeling, driver behavior associated with advanced traveler information messages, transportation planning for connected and autonomous vehicles, performance assessment for alternative intersections and interchanges, and geometric design. Dr. Yang serves as an editorial board member for the Journal of "Case Studies on Transport Policies", "International Journal of Transportation Science and Technology", and "Transportation Safety and Environmental (English Edition)" and is a young academic editor of the "Journal of Traffic and Transportation Engineering (English Edition)".

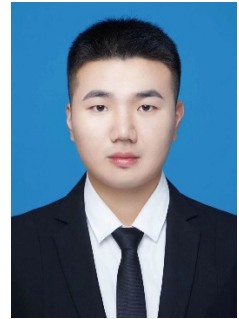

Chuang Zhu (Graduate Student Member, IEEE) received the B.S. degree in transportation from Shijiazhuang Tiedao University in 2020. He is currently working toward the Ph.D. degree in School of Traffic and Transportation, Beijing Jiaotong University, Beijing, China. His research interests include public transport optimization and urban transport demand management.

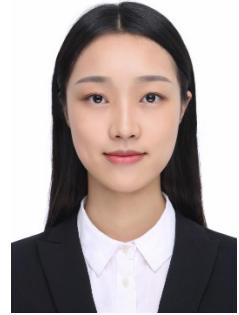

Siyuan Chen graduated from Beijing Jiaotong University in 2020 with a Bachelor's degree in Traffic Engineering and graduated from Beijing Jiaotong University in 2023 with a Master's degree in Traffic Planning and Management. She is currently working in Zhejiang Communications Investment Group. Her research interests include traffic planning and ridesharing route optimization.